\documentclass[12pt]{amsart}
\usepackage{fullpage,url}
\usepackage[all]{xy} 

\DeclareFontEncoding{OT2}{}{} 

\usepackage{color}

\newcommand{\defi}[1]{\textsf{#1}} 

\newcommand{\Aff}{{\mathbb A}}

\newcommand{\F}{{\mathbb F}}
\newcommand{\G}{{\mathbb G}}

\newcommand{\Q}{{\mathbb Q}}
\newcommand{\R}{{\mathbb R}}
\newcommand{\Z}{{\mathbb Z}}

\newcommand{\frakA}{{\mathfrak A}}
\newcommand{\frakB}{{\mathfrak B}}
\newcommand{\frakH}{{\mathfrak H}}

\newcommand{\mm}{{\mathfrak m}}

\newcommand{\calA}{{\mathcal A}}

\newcommand{\OO}{{\mathcal O}}

\DeclareMathOperator{\good}{good}
\DeclareMathOperator{\Hilb}{Hilb}

\DeclareMathOperator{\Tr}{Tr}

\DeclareMathOperator{\Aut}{Aut}

\DeclareMathOperator{\Gr}{Gr}

\DeclareMathOperator{\Sym}{Sym}

\DeclareMathOperator{\Spec}{Spec}
\DeclareMathOperator{\Isom}{Isom}
\DeclareMathOperator{\SPEC}{\bf Spec}

\newcommand{\Schemes}{\operatorname{\bf Schemes}}
\newcommand{\Sets}{\operatorname{\bf Sets}}

\newcommand{\op}{{\operatorname{op}}}

\newcommand{\et}{{\operatorname{et}}}
\newcommand{\spl}{{\operatorname{split}}}

\newcommand{\tH}{{\operatorname{th}}}

\newcommand{\GL}{\operatorname{GL}}

\newcommand{\injects}{\hookrightarrow}
\newcommand{\isom}{\simeq}

\newcommand{\intersect}{\cap} 
\newcommand{\tensor}{\otimes}
\newcommand{\directsum}{\oplus} 

\numberwithin{equation}{section}

\newtheorem{theorem}[equation]{Theorem}
\newtheorem{lemma}[equation]{Lemma}
\newtheorem{corollary}[equation]{Corollary}
\newtheorem{proposition}[equation]{Proposition}

\theoremstyle{definition}
\newtheorem{definition}[equation]{Definition}
\newtheorem{warning}[equation]{Warning}
\newtheorem{conjecture}[equation]{Conjecture}
\newtheorem{question}[equation]{Question}
\newtheorem{example}[equation]{Example}
\newtheorem{examples}[equation]{Examples}

\theoremstyle{remark}
\newtheorem{remark}[equation]{Remark}

\usepackage[
        backref,
        pdfauthor={Bjorn Poonen}, 
]{hyperref}
\usepackage[alphabetic,backrefs,lite]{amsrefs} 

\begin{document}

\title[Moduli space of algebras]{The moduli space of commutative algebras of finite rank}
\subjclass[2000]{Primary 14D20; Secondary 13E10, 13M99, 14C05}
\author{Bjorn Poonen}
\address{Department of Mathematics, University of California, 
        Berkeley, CA 94720-3840, USA}
\email{poonen@math.berkeley.edu}
\urladdr{http://math.berkeley.edu/\~{}poonen/}
\date{March 8, 2007}

\begin{abstract}
The moduli space of rank-$n$ commutative algebras 
equipped with an ordered basis
is an affine scheme $\frakB_n$ of finite type over $\Z$,
with geometrically connected fibers.
It is smooth if and only if $n \le 3$.
It is reducible if $n \ge 8$ 
(and the converse holds, at least if we remove the fibers above $2$ and $3$).
The relative dimension of $\frakB_n$ is $\frac{2}{27} n^3 + O(n^{8/3})$.
The subscheme parameterizing \'etale algebras
is isomorphic to $\GL_n/S_n$,
which is of dimension only $n^2$.
For $n \ge 8$, there exist algebras are not limits of \'etale algebras.
The dimension calculations lead also to new asymptotic formulas 
for the number of commutative rings of order $p^n$
and the dimension of the Hilbert scheme of $n$ points
in $d$-space for $d \ge n/2$.
\end{abstract}

\maketitle

\section{The moduli space of based algebras}
\label{S:representable}

All rings and algebras are assumed to be 
commutative, associative, and with $1$,
unless otherwise specified.

Let $n \in \Z_{\ge 0}$.
For every ring $k$, we would like to parameterize the $k$-algebras
that are locally free of rank $n$ as a $k$-module.
But for $n \ge 2$ 
this moduli problem turns out to be not representable by a scheme;
the difficulty is related to presence of automorphisms.
To rigidify, we equip $k$-algebras with extra structure.
One option would be to consider $k$-algebras equipped with a ordered set
of $d$ generators as a $k$-algebra: these are parameterized by
the Hilbert scheme $\Hilb^n(\Aff^d)$ of $n$ points in $\Aff^d$.
Another option, the one explored in this paper, is
to consider $k$-algebras equipped with an ordered basis.

Define a functor $\frakB_n \colon \Schemes^\op \to \Sets$
as follows:
\begin{enumerate}
\item
If $S$ is a scheme, an element of $\frakB_n(S)$ is a pair $(\calA,\phi)$
(strictly speaking, an isomorphism class of pairs)
where $\calA$ is an $\OO_S$-algebra
and $\phi\colon \calA \to \OO_S^{\directsum n}$
is an $\OO_S$-module isomorphism.
Two pairs $(\calA,\phi)$ and $(\calA',\phi')$ are considered the same
if there is an isomorphism $\calA \to \calA'$ that with $\phi$ and $\phi'$
makes a commuting triangle.
\item
If $f \colon T \to S$ is a morphism of schemes,
then $\frakB_n(f)\colon \frakB_n(S) \to \frakB_n(T)$
maps $(\calA,\phi)$ to $(f^* \calA,f^*\phi)$.
\end{enumerate}

\begin{proposition}
The functor $\frakB_n$ 
is representable by an affine scheme of finite type over~$\Z$.
\end{proposition}

\begin{proof}
Let $e_1,\ldots,e_n$ be the standard $\OO_S$-basis of $\OO_S^{\directsum n}$.
Given $(\calA,\phi) \in \frakB_n(S)$,
we may use $\phi$ to transport the multiplication on $\calA$
to a multiplication on $\OO_S^{\directsum n}$,
which is described by giving the $n^3$ constants $c_{ij}^\ell \in H^0(S,\OO_S)$
such that $e_i e_j = \sum_{\ell=1}^n c_{ij}^\ell e_\ell$;
we also get $\phi(1) = \sum d_i e_i$
for some $d_i \in H^0(S,\OO_S)$.
Conversely, given elements $c_{ij}^\ell,d_i \in H^0(S,\OO_S)$
satisfying the polynomial conditions 
that say that the resulting multiplication law
is commutative and associative 
and has $\sum d_i e_i$ as a multiplicative identity,
we can recover the algebra uniquely.
Thus $\frakB_n$ is representable by the closed subscheme of $\Aff^{n^3+n}_\Z$
cut out by these polynomial conditions
on {\em indeterminates} $c_{ij}^\ell$ and $d_i$.
\end{proof}

{}From now on, $\frakB_n$ denotes the representing scheme.
Call $\frakB_n$ the \defi{moduli space of based rank-$n$ algebras}.

\begin{remark}
Various analogues of $\frakB_n$ have been studied in the literature
before, at least over a base field.
The paper \cite{Flanigan1968} defines the moduli space for based
associative (but not necessarily commutative or unital) algebras, 
defines a $\GL_n$-action as we do below,
and studies $\GL_n$-orbits in its irreducible components.
The moduli space of based associative algebras with~$1$,
along with the $\GL_n$-action, is studied 
in \cites{Gabriel1974,Mazzola1979,Mazzola1982,LeBruyn-Reichstein1999}.
In particular, \cite{Mazzola1979} contains a detailed study
of its number of irreducible components
and computes the number explicitly for $n=5$ (it is $10$).
The paper \cite{Mazzola1980}
studies the moduli space of based nilpotent 
commutative associative rank-$n$ algebras over an algebraically closed field
of characteristic not $2$ or $3$, especially for $n \le 6$;
it follows from this work that $\frakB_n$ has geometrically irreducible
fibers, at least over $\Z[1/6]$, for all $n \le 7$.
Asymptotic formulas for the dimension of $\frakB_n$ 
and of the moduli space of based associative algebras with $1$
are mentioned in \cite{Neretin1987}, whose primary goal
was an asymptotic formula for the dimension of 
the moduli space of $n$-dimensional Lie algebras.
Finally, the moduli space $\frakB_n$ is implicit in some of the work
of Bhargava \cites{Bhargava2004I,Bhargava2004II,Bhargava2004III,Bhargava2005}
on parameterizing and enumerating algebras of rank $\le 5$,
especially over $\Z$.
\end{remark}

\begin{definition}
Given $S$, define two points of $\frakB(S)$ as follows:
\begin{enumerate}
\item
$\tilde{\calA}_{\spl}$ is $\calA_{\spl}:=\OO_S \times \cdots \times \OO_S$
equipped with the obvious basis.
\item
$\bullet$ is 
$\calA_\bullet:=\frac{\OO_S[x_1,\ldots,x_{n-1}]}{(x_1,\ldots,x_{n-1})^2}$
equipped with the basis $1,x_1,\ldots,x_{n-1}$.
\end{enumerate}
\end{definition}

For some calculations it will be simpler to work with based algebras
where the first basis vector equals the multiplicative identity:
\begin{definition}
For $n \ge 1$, let $\frakB_n^1$ be the scheme
parameterizing based algebras such that the first basis vector equals $1$.
It is the closed subscheme of $\frakB_n$ cut out 
by the equations $d_1=1$ and $d_i=0$ for $i \ge 2$.
For $n=0$, we use the convention $\frakB_0^1 := \frakB_0$.
\end{definition}

\section{Action of the general linear group}
\label{S:GL_n}

Let $\GL_n$ be the general linear group scheme over $\Z$.
Each $M \in \GL_n(S) = \Aut_{\OO_S}(\OO_S^{\directsum n})$
acts on $\frakB_n(S)$ by sending $(\calA,\phi)$ to $(\calA,M \circ \phi)$.
This defines a left action of $\GL_n$ on $\frakB_n$.

Given two elements 
$(\calA,\phi)$ and $(\calA',\phi')$ of $\frakB_n(S)$,
the set $\Isom(\calA,\calA')$
of $\OO_S$-algebra isomorphisms $\alpha\colon \calA \to \calA'$
is in bijection with the set of $M \in \GL_n(S)$
mapping $(\calA,\phi)$ to $(\calA',\phi')$.
Namely, $\alpha$ corresponds to the $M \in \GL_n(S)$ making
\[
\xymatrix{
\calA \ar[r]^{\phi} \ar[d]_{\alpha} & \OO_S^{\directsum n} \ar[d]^M \\
\calA' \ar[r]^{\phi'} & \OO_S^{\directsum n} \\
}
\]
commute.

It follows that the set of isomorphism classes 
of free rank-$n$ $\OO_S$-algebras 
is in bijection with the quotient set $\GL_n(S)\backslash \frakB_n(S)$.
Also, given $(\calA,\phi) \in \frakB_n(S)$,
the group of $\OO_S$-algebra automorphisms $\Aut(\calA)$
is identified with the stabilizer in $\GL_n(S)$ of 
$(\calA,\phi) \in \frakB_n(S)$.

\section{Comparison of $\frakB_n$ and $\frakB_n^1$}
\label{S:contracted product}

For $n \ge 1$,
let $H$ be the subgroup scheme of $\GL_n$ stabilizing $(1,0,\ldots,0)$.
It consists of invertible matrices whose first column equals 
$(1,0,\ldots,0)$.
If $n=0$, we use the convention that $H$, like $\GL_n$, is 
the trivial group scheme.
In any case, the relative dimension of $H$ is $n^2-n$.

The left action of $\GL_n$ on $\frakB_n$ restricts
to a left action of $H$ on $\frakB_n^1$.
Also, there is a right action of $H$ on $\GL_n$ given by multiplication
on the right.
Combining these gives a left $H$-action on 
$\GL_n \times \frakB_n^1$
in which $N \in H(S)$ maps
$(M,\tilde{\calA}) \in (\GL_n \times \frakB_n^1)(S)$
to
$(M N^{-1},N \cdot \tilde{\calA})$.
The \defi{contracted product} $\GL_n \overset{H}\wedge \frakB_n^1$
is the quotient of $\GL_n \times \frakB_n^1$ by this (free) $H$-action.

\begin{proposition}
\label{P:contracted product}
The scheme $\frakB_n$ is isomorphic to $\GL_n \overset{H}\wedge \frakB_n^1$.
\end{proposition}

\begin{proof}
Restricting the action $\GL_n \times \frakB_n \to \frakB_n$
yields a map
\[
        \GL_n \times \frakB_n^1 \to \frakB_n.
\]
which is invariant for the $H$-action on $\GL_n \times \frakB_n^1$.
The fiber above $\tilde{\calA} \in \frakB_n(S)$
in $\GL_n \times \frakB_n^1$
is a Zariski locally trivial torsor under $H$,
trivialized by any open covering of $S$ 
that splits the $\OO_S$-module injection 
$\OO_S \to \calA$ given by the algebra structure.
\end{proof}

\begin{definition}
Define the \defi{moduli stack of locally free rank-$n$ algebras}
as the Artin stack
\[
        \frakA_n := [\GL_n \backslash \frakB_n]  
        \isom [H \backslash \frakB_n^1].
\]
(Here $\backslash$ denotes 
``stack quotient by the group acting on the left'',
not set difference!)
\end{definition}

\section{Comparison with Hilbert schemes}
\label{S:Hilbert}

For any $n,d \in \Z_{\ge 0}$,
let $\Hilb_n(\Aff^d)$ be the \defi{Hilbert scheme of $n$ points in $\Aff^d$}:
it is the $\Z$-scheme whose $S$-points correspond
(functorially in $S$) to closed subschemes $T$ of $\Aff^d_S$ that are flat
over $S$ and whose geometric fibers have Hilbert polynomial
equal to the constant polynomial $n$.
Let $\frakH_n(\Aff^d)$ 
be the $\Z$-scheme that parameterizes closed subschemes $T$ as above
equipped with an ordered $\OO_S$-basis for $\OO_T$.
Then $\frakH_n(\Aff^d)$ is a torsor under $\GL_n$ over $\Hilb_n(\Aff^d)$.
There is a morphism $\frakH_n(\Aff^d) \to \frakB_n$
that maps $T \subseteq \Aff^d$ with an $\OO_S$-basis
to $\OO_T$ with the $\OO_S$-basis (forgetting the embedding into $\Aff^d$).
This is $\GL_n$-equivariant, so we obtain a cartesian diagram
\[
\xymatrix{
\frakH_n(\Aff^d) \ar[r] \ar[d] & \frakB_n \ar[d] \\
\Hilb_n(\Aff^d) \ar[r] & \frakA_n. \\
}
\]

We may view $\frakH_n(\Aff^d)$ as the moduli space parameterizing rank-$n$
$\OO_S$-algebras $\OO_T$ equipped with both a basis
and a sequence $(x_1,\ldots,x_d)$ of global sections
that generate $\OO_T$ as an $\OO_S$-algebra
(to give these global sections is 
to give a closed immersion $T \injects \Aff^d_S$).
On the other hand, $\frakB_n \times \Aff^{nd}$ 
is the moduli space parameterizing rank-$n$ $\OO_S$-algebras $\OO_T$
equipped with both a basis and
an {\em arbitrary} sequence $(x_1,\ldots,x_d)$ of global sections.
(The $nd$ coordinates on $\Aff^{nd}$ specify the coordinates of
the $x_i$ with respect to the basis.)

\begin{proposition}
\label{P:Hilbert open}
The morphism 
$\frakH_n(\Aff^d) \to \frakB_n \times \Aff^{nd}$ is an open immersion.
\end{proposition}

\begin{proof}
Over a field $k$, the vector space spanned by the monomials of total
degree $\le e$ in $x_1,\ldots,x_d$ increases strictly in dimension
with $e$ until it stabilizes, so it must stabilize at or before $e=n-1$.
Thus $x_1,\ldots,x_d$ generate the rank-$n$ $k$-algebra $A$
if and only if the monomials of total degree $<n$ span $A$.
This is a condition on the rank of a matrix whose entries
are polynomials in the coordinates of the $x_i$
with respect to a fixed basis of $A$;
the condition is expressible as the nonvanishing of certain minors.
Those minors, viewed as functions on $\frakB_n \times \Aff^{nd}$,
cut out the complement of $\frakH_n(\Aff^d)$ in $\frakB_n \times \Aff^{nd}$.
\end{proof}

\begin{corollary}
\label{C:open locus in Bn}
The image of $\frakH_n(\Aff^d) \to \frakB_n$ is open in $\frakB_n$.
\end{corollary}

\begin{proposition}
\label{P:Hilbert section}
If $d \ge n$, then 
$\frakH_n(\Aff^d) \to \frakB_n$ admits a section.
The map $\frakH_n(\Aff^d) \to \frakB_n$
is surjective if and only if $d \ge n-1$.
\end{proposition}

\begin{proof}
If $d \ge n$, define $\frakB_n \to \frakH_n(\Aff^d)$
by mapping an $\OO_S$-algebra $\calA$ with basis $e_1,\ldots,e_n$,
to the same $\OO_S$-algebra with the same basis
and with algebra generators $(x_1,\ldots,x_d):=(e_1,\ldots,e_n,0,0,\ldots,0)$;
this is a section.

Given a based $\OO_S$-algebra,
Zariski locally on $S$ we may replace one of the basis elements with $1$
and still have a basis, and then the other basis elements are algebra
generators.
This proves surjectivity for $d \ge n-1$.

We will not have surjectivity for $d<n-1$, 
since $\calA_\bullet$ cannot be generated 
as an algebra by fewer than $n-1$ elements.
\end{proof}

Let $k$ be a field.
For any $\Z$-scheme $V$, let $V_k = V \underset{\Spec \Z} \times \Spec k$.
For instance, $\frakB_{n,k}$ denotes the $k$-variety (not
necessarily irreducible) obtained from $\frakB_n$.

\begin{corollary}
\label{C:Hilbert scheme dimension}
We have
\[
	\dim \Hilb_n(\Aff^d)_k \le \dim \frakB_{n,k} - n^2 + nd.
\]
If $d \ge n-1$, then equality holds.
\end{corollary}

\begin{proof}
By Propositions \ref{P:Hilbert open},
\[
	\dim \Hilb_n(\Aff^d)_k = \dim \frakH_n(\Aff^d)_k - \dim \GL_{n,k}
	\le \dim \frakB_{n,k} + \dim \Aff^{nd}_k - \dim \GL_{n,k}.
\]
By~\ref{P:Hilbert section}, equality holds if $d \ge n-1$.
\end{proof}

\begin{remark}
If we are interested in only the case $n=d$, then the connection
between Hilbert schemes and $\frakB_n$ is even more direct.
There is an open subscheme $H_n$ of $\Hilb_n(\Aff^n)$
consisting of the points for which the corresponding quotient
of the polynomial algebra is spanned as a vector space by the 
coordinate functions on $\Aff^n$.
Also, for any rank-$n$ algebra with a fixed basis, the locus of bases 
is contained in the locus of $n$-tuples of algebra generators,
which is contained in the affine space of all $n$-tuples of elements,
and the first is dense in the third, so the first is dense in the second,
which implies that the open subscheme $H_n$ of $\Hilb_n(\Aff^n)$
is Zariski dense.
Now we have an isomorphism $\frakB_n \isom H_n$
taking a based algebra to the same algebra equipped with the $n$ algebra
generators given by the basis.
This isomorphism may also be viewed as the composition of the section 
$\frakB_n \to \frakH_n(\Aff^n)$
of Proposition~\ref{P:Hilbert section} 
with the morphism $\frakH_n(\Aff^n) \to \Hilb_n(\Aff^n)$.
The fact that $\frakB_n$ is isomorphic to a dense open subscheme
of $\Hilb_n(\Aff^n)$ explains why many of the properties
of $\frakB_n$ we prove in later sections 
reflect the corresponding properties of Hilbert schemes.
\end{remark}

\begin{remark}
\label{R:smooth up to 3}
By \cite{Fogarty1968} (see especially Corollary~2.10 and the proof
of Theorem~2.4 there),
$\Hilb_n(\Aff^2)$ is smooth over $\Z$ with geometrically irreducible fibers.
It follows that the same is true of $\frakH_n(\Aff^2)$,
and (by Proposition~\ref{P:Hilbert open})
of the image of $\frakH_n(\Aff^2) \to \frakB_n$
viewed as an open subscheme of $\frakB_n$
(the locus parameterizing based algebras that can be locally generated
as an algebra by $\le 2$ elements).
If $n \le 3$, Proposition~\ref{P:Hilbert section} 
implies that this locus is the whole scheme $\frakB_n$,
which is then smooth over $\Z$ with geometrically irreducible fibers.
\end{remark}

\begin{remark}
Fix a characteristic $p \ge 0$ and a nonnegative integer $d$.
The value of $\dim \Hilb_n(X)$ is the same for every 
nonempty smooth variety $X$ of dimension $d$ over a field of
characteristic~$p$.
\end{remark}

\section{The moduli spaces for $n \le 3$}
\label{S:small n}

\begin{proposition}
\label{P:small n}
For $n \le 3$, the isomorphism type of the $\Z$-scheme $\frakB_n^1$ 
is given by the following table:
\begin{center}
\begin{tabular}{c|cccc}
$n$ & $0$ & $1$ & $2$ & $3$ \\ \hline
$\frakB_n^1$ & $\Aff^0$ & $\Aff^0$ & $\Aff^2$ & $\Aff^6$
\end{tabular}
\end{center}
\end{proposition}

\begin{proof}
On $\frakB_n^1 \subseteq \frakB_n \subseteq \Aff^{n^3+n}$ 
the values of the coordinates $d_i$ are specified,
and the values of $c_{ij}^\ell$ when $i$ or $j$ is $1$
are forced by the defining polynomials.
Thus we may view $\frakB_n^1$ as a closed subscheme in
the affine space $\Aff^{(n-1)^2 n}$
whose coordinates are the indeterminates $c_{ij}^\ell$
for $2 \le i,j \le n$ and $1 \le \ell \le n$.

When $n \le 2$, the commutativity and associativity conditions
on these $c_{ij}^\ell$ are vacuous,
so $\frakB_n^1 \isom \Aff^{(n-1)^2 n}$.

Finally, suppose $n=3$.
Call a basis for a rank-$3$ algebra \defi{good}
if it has the form $1,\alpha,\beta$ with $\alpha \beta \in \OO_S \cdot 1$.
We paraphrase \cite{Delone-Faddeev1940}*{\S15} 
(see also \cite{Gan-Gross-Savin2002}*{\S4}),
where it was shown that rank-$3$ algebras (over $\Z$)
equipped with a good basis are in bijection with
binary cubic forms.
Let $\frakB_n^{1,\good}$ be the closed subscheme of $\frakB_n^1$
parameterizing rank-$3$ algebras equipped with a good basis.
There is a left action of $\G_a^2$ on $\frakB_n^1$ in which
$(a,b) \in \G_a^2(S)$ maps an algebra $\calA$ equipped
with basis $1,\alpha,\beta$
to $\calA$ equipped with basis $1,\alpha+a,\beta+b$.
This action restricts to an isomorphism
\[
        \G_a^2 \times \frakB_n^{1,\good} \to \frakB_n^1.
\]
Using the conditions coming from commutativity and associativity,
one finds that the multiplication in any based algebra in $\frakB_n^{1,\good}$
with basis $1,\alpha,\beta$ has the form
\begin{align*}
        \alpha^2 &= -ac + b \alpha - a \beta \\
        \beta^2 &= -bd + d \alpha - c \beta \\
        \alpha \beta &= -ad
\end{align*}
for some $a,b,c,d$,
and that conversely any $a,b,c,d$ yield a based algebra 
in $\frakB_n^{1,\good}$.
Thus $\frakB_n^{1,\good} \isom \Aff^4$,
and $\frakB_n^1 \isom \Aff^2 \times \Aff^4 \isom \Aff^6$.
\end{proof}

\begin{remark}
In Section~\ref{S:smoothness} we will find that for $n \ge 4$,
the scheme $\frakB_n^1$ is not smooth over $\Spec \Z$,
so in particular it is not isomorphic to affine space.
\end{remark}

\section{The \'etale locus}
\label{S:etale}

Given $\tilde{\calA}:=(\calA,\phi) \in \frakB_n(S)$,
we may define an $\OO_S$-linear trace map $\Tr\colon \calA \to \OO_S$
in the usual way.
Let $b_i =\phi^{-1}(e_i) \in H^0(S,\calA)$;
these form an $\OO_S$-basis of $\calA$.
Define a regular function $\Delta \colon \frakB_n \to \Aff^1$ by setting
\[
        \Delta(\tilde{\calA}) := \det \left(\Tr(b_i b_j)\right) 
        \in H^0(S,\OO_S).
\]
Acting on $\tilde{\calA}$ by $M \in \GL_n(S)$
multiplies the matrix in the definition of $\Delta(\tilde{\calA})$
by $(M^{-1})^t$ on the left and $M^{-1}$ on the right, so
\[
        \Delta(M \cdot \tilde{\calA}) = (\det M)^{-2} \Delta(\tilde{\calA}).
\]
In particular, the zero locus $\{\Delta=0\}$ and its complement $\frakB_n^\et$
in $\frakB_n$ are $\GL_n$-invariant.
Call the open affine subscheme $\frakB_n^\et$ of $\frakB_n$ 
the \defi{\'etale locus}.
The following is standard:

\begin{proposition}
\label{P:etale locus}
The following are equivalent for 
$\tilde{\calA} := (\calA,\phi) \in \frakB_n(S)$:
\begin{enumerate}
\item[(i)]
$\tilde{\calA} \in \frakB_n^\et(S)$.
\item[(ii)]
$\Delta(\tilde{\calA}) \in H^0(S,\OO_S^\times)$.
\item[(iii)]
The morphism $\SPEC \calA \to S$ is \'etale.
\item[(iv)]
There exists a surjective \'etale base extension $f\colon T \to S$
and an isomorphism $f^* \calA \isom f^* \calA_\spl$.
\end{enumerate}
\end{proposition}

\begin{examples}\hfill
\begin{enumerate}
\item
We have $\Delta(\tilde{\calA}_\spl)=1$, 
so $\tilde{\calA}_\spl \in \frakB_n^\et(S)$.
\item
If $n \ge 2$, $\Delta(\bullet)=0$,
so $\bullet \notin \frakB_n^\et(S)$.
\end{enumerate}
\end{examples}

Let $S_n$ denote the constant group scheme over $\Z$ 
corresponding to the symmetric group on $n$ letters.
Embed $S_n$ in $\GL_n$ as the subgroup of permutation matrices.

\begin{theorem}
\label{T:homogeneous space}
There is a $\GL_n$-equivariant isomorphism from the homogeneous
space $\GL_n/S_n$ to $\frakB_n^\et$.
\end{theorem}

\begin{proof}
We have a morphism $\GL_n \to \frakB_n^\et$
sending $M$ to $M \cdot \tilde{\calA}_\spl$.
The equivalence of (i) and (iv) in Proposition~\ref{P:etale locus}
implies that $\GL_n \to \frakB_n^\et$ is surjective.

It remains to show that the subgroup scheme 
of $\GL_n$ stabilizing $\tilde{\calA}_\spl$ equals $S_n$.
Equivalently, we must show that for any connected scheme $S$,
the automorphism group of the $\OO_S$-algebra $\prod_{i=1}^n \OO_S$
equals $S_n$ (acting by permuting coordinates).
This holds, since every automorphism of the $S$-scheme
$\coprod_{i=1}^n S$ induces a permutation of the connected components,
and each map between components is the identity $S \to S$
by virtue of being an $S$-morphism.
\end{proof}

\begin{corollary}
\label{C:et is irreducible}
The scheme $\frakB_n^\et$ is irreducible.
\end{corollary}

\begin{corollary}
The Zariski closure $\overline{\frakB_n^\et}$ of
$\frakB_n^\et$ in $\frakB_n$ is an irreducible component of $\frakB_n$.
\end{corollary}

\begin{corollary}
The affine variety $\frakB_{n,k}^\et$ is irreducible and of dimension $n^2$.
\end{corollary}

\begin{corollary}
\label{C:irreducible component}
The Zariski closure $\overline{\frakB_{n,k}^\et}$ 
of $\frakB_{n,k}^\et$ in $\frakB_{n,k}$
is an irreducible component of $\frakB_{n,k}$ of dimension $n^2$.
\end{corollary}

\begin{warning}
Conceivably, the base extension
$\left( \overline{\frakB_n^\et} \right)_k$
could be strictly larger than $\overline{\frakB_{n,k}^\et}$.
Hence we cannot conclude that 
$\left( \overline{\frakB_n^\et} \right)_k$
is irreducible.
\end{warning}

\begin{remark}
We could similarly define $\frakB_n^{1,\et}$
and prove that $\frakB_n^{1,\et} \isom H/S_{n-1}$,
where $S_{n-1} := S_n \intersect H \subseteq \GL_n$.
The analogues of 
Corollaries \ref{C:et is irreducible} to~\ref{C:irreducible component}
for $\frakB_n^{1,\et}$ follow.
In particular, $\dim \frakB_{n,k}^{1,\et} = n^2-n$.
\end{remark}

\begin{remark}
The number of isomorphism types of rank-$n$ algebras
over an algebraically closed field $k$ is finite
if and only if $n \le 6$ \cite{Poonen2007-dimension6-preprint}.
It follows that $\dim \frakB_{n,k}=n^2$
and $\dim \frakB_{n,k}=n(n-1)$
for $n \le 6$, even if $k$ is not algebraically closed.
Actually, these dimension formulas hold also for $n=7$,
at least if $k$ is of characteristic not $2$ or $3$,
because $\frakB_{n,k}$ is irreducible by \cite{Mazzola1980}*{Corollary~4},
and hence equals $\overline{\frakB_n^\et}$.
\end{remark}

\begin{remark}
One might expect at first that all $k$-algebras are degenerations
of \'etale algebras; i.e., that $\frakB_{n,k}^\et$ is Zariski dense
in $\frakB_{n,k}$.
In Section~\ref{S:dimension lower bound} 
we will disprove this for large $n$
by proving that $\dim \frakB_{n,k} \ge \frac{2}{27} n^3$.
In fact, for every $n \ge 8$, the variety $\frakB_{n,k}$ 
has irreducible components other than $\overline{\frakB_n^\et}$:
see Proposition~\ref{P:8}.
\end{remark}

\section{Connectedness}
\label{S:connectness}

\begin{proposition}
\label{P:cone}
Suppose $n \ge 1$.
Then $\frakB_n^1$ is the affine cone over a closed subscheme
in a weighted projective space.
The vertex of the cone is $\bullet$.
\end{proposition}

\begin{proof}
As in the proof of Proposition~\ref{P:small n},
view $\frakB_n^1$ as a closed subscheme of
the $\Aff^{(n-1)^2 n}$ with coordinates 
$c_{ij}^\ell$ with $2 \le i, j \le n$ and $1 \le \ell \le n$.
Let the weight of $c_{ij}^\ell$ be $2$ if $\ell=1$, and $1$ otherwise.
Then the equations expressing commutativity and associativity 
are homogeneous.
The origin in $\Aff^{(n-1)^2 n}$
corresponds to the multiplication table for $\bullet$.
\end{proof}

\begin{corollary}
The point $\bullet \in \frakB_n^1(k)$ 
belongs to every irreducible component of $\frakB_{n,k}^1$.
\end{corollary}

\begin{remark}
The observation that every algebra in $\frakB_{n,k}^1$
can be connected to $\bullet$ 
is due to Manjul Bhargava.
\end{remark}

\begin{corollary}
The point $\bullet$
belongs to every irreducible component of $\frakB_{n,k}$.
\end{corollary}

\begin{proof}
The inclusion $\frakB_n^1 \to \frakB_n$ factors as
\begin{equation}
\label{E:factorization}
\frakB_n^1 \longrightarrow \GL_n \times \frakB_n^1 \longrightarrow \frakB_n
\end{equation}
in which the first map sends $\tilde{\calA}$
to $(1_n,\tilde{\calA})$,
and the second map is a torsor under $H$,
by Proposition~\ref{P:contracted product}.
Since $\GL_n$ and $H$ are irreducible,
the base extension of \eqref{E:factorization} to $k$
induces bijections on irreducible components.
\end{proof}

As a corollary, we obtain

\begin{theorem}
\label{T:connected}
The varieties $\frakB_{n,k}^1$ and $\frakB_{n,k}$ are connected.
\end{theorem}

\section{Smoothness}
\label{S:smoothness}

\begin{proposition}
For $n \ne 2$,
the tangent space $T_{\bullet}(\frakB_{n,k}^1)$ 
has dimension $n(n-1)^2/2$.
\end{proposition}

\begin{proof}
By Proposition~\ref{P:small n},
we may assume $n \ge 3$.
A deformation of $\bullet$ to a based algebra
in $\frakB_{n,k}^1(k[\epsilon]/(\epsilon^2))$
is given in terms of a basis $e_1,\ldots,e_n$ (with $e_1=1$)
by the $\frac{n(n-1)}{2} \cdot n$ multiplication constants $c_{ij}^\ell$
for $2 \le i \le j \le n$ and $1 \le \ell \le n$,
in which $c_{ij}^\ell = \gamma_{ij}^\ell \epsilon$
with $\gamma_{ij}^\ell \in k$.
In assuming $i \le j$, we have already imposed commutativity,
so it remains to examine the restrictions on the
$\gamma_{ij}^\ell$ imposed by associativity.

The condition $(e_2 e_2) e_3 = e_2 (e_2 e_3)$
implies $\gamma_{22}^1 = 0$ and $\gamma_{23}^1 = 0$.
Similarly $\gamma_{ij}^1=0$ for all $i,j$.
No other conditions are imposed, so the dimension
equals the number of $(i,j,\ell)$ satisfying
$2 \le i \le j \le n$ and $2 \le \ell \le n$.
This is $\frac{n(n-1)}{2} \cdot (n-1)$.
\end{proof}

\begin{corollary}
For $n \ne 2$,
the tangent space $T_{\bullet}(\frakB_{n,k})$
has dimension $n(n-1)^2/2 + n$.
\end{corollary}

\begin{proof}
Using \eqref{E:factorization}, we find 
$\dim T_{\bullet}(\frakB_{n,k}) 
        = \dim T_{\bullet}(\frakB_{n,k}^1) + \dim \GL_n - \dim H.$
\end{proof}

\begin{corollary}
\label{C:singular point}
For $n \ge 4$, the point $\bullet$
is singular on both $\frakB_{n,k}^1$ and $\frakB_{n,k}$.
\end{corollary}

\begin{proof}
The irreducible component $\overline{\frakB_{n,k}^{1,\et}}$
of $\frakB_{n,k}^1$ contains $\bullet$
and its dimension $n(n-1)$ is less than 
$\dim T_{\bullet}(\frakB_{n,k}^1) = n(n-1)^2/2$
if $n \ge 4$.
\end{proof}

\begin{proposition}
For $n \in \Z_{\ge 0}$, the following are equivalent:
\begin{enumerate}
\item $n \le 3$.
\item $\frakB_{n,k}^1$ is smooth.
\item $\frakB_{n,k}$ is smooth.
\end{enumerate}
\end{proposition}

\begin{proof}
By Proposition~\ref{P:contracted product}, $\frakB_{n,k}^1$ is smooth
if and only if $\frakB_{n,k}$ is.
For $n\le 3$ use Proposition~\ref{P:small n} 
or Remark~\ref{R:smooth up to 3}.
For $n \ge 4$ use Corollary~\ref{C:singular point}.
\end{proof}

\begin{example}
Let $R$ be a discrete valuation ring with uniformizer $\pi$.
Let $A$ be the based $R/\pi^2 R$-algebra with basis $1,x,y,z$
satisfying $x^2=\pi x$, $y^2=\pi y$, $z^2=\pi z$, $xy=\pi z$, 
$yz=\pi x$, $zx=\pi y$.
Then $A$ does not lift to a free $R/\pi^3 R$-algebra of rank $4$,
let alone a free $R$-algebra of rank $4$:
the associative relation $x^2 y = x(xy)$ fails for any lift to $R/\pi^3 R$.
This violation of the infinitesimal lifting criterion for smoothness
shows again that $\frakB_4$ is not smooth.
\end{example}

\section{Dimension: lower bound}
\label{S:dimension lower bound}

\begin{lemma}
\label{L:Zdr}
Suppose $n,d,r \in \Z_{\ge 0}$ satisfy $r \le d(d+1)/2$ and $n=1+d+r$.
Then
\[
\dim \frakB_{n,k} \ge r \left( \frac{d(d+1)}{2} - r \right) + n^2 - (d^2+dr).
\]
\end{lemma}

\begin{proof}
Let $\mm=(x_1,\ldots,x_d) \subseteq k[x_1,\ldots,x_d]$.
Then $\dim \mm^2/\mm^3 = d(d+1)/2$.
Let $Z_{d,r}$ be the variety parameterizing $(V,\phi)$
where $V$ is a codimension-$r$ subspace of $\mm^2/\mm^3$,
and $\phi$ is a $k$-linear isomorphism 
$\frac{k[x_1,\ldots,x_d]}{\mm^3}/V \to k^n$.
The forgetful map $(V,\phi) \mapsto V$ exhibits $Z_{d,r}$
as a $\GL_n$-torsor over the Grassmannian $\Gr(r,d(d+1)/2)$,
so
\[
        \dim Z_{d,r} = \dim \Gr(r,d(d+1)/2) + \dim \GL_n 
        = r \left( \frac{d(d+1)}{2} - r \right) + n^2.
\]

We next compute the dimension of the fibers of 
the morphism $\pi\colon Z_{d,r} \to \frakB_{n,k}$
defined by $(V,\phi) \mapsto (A_V,\phi)$
where $A_V:=\frac{k[x_1,\ldots,x_d]}{\mm^3}/V$.
For fixed $V$, to give a $k$-algebra isomorphism from $A_V$ to some $A_{V'}$
is to give $(\ell,\eta_1,\ldots,\eta_d)$
where $\ell$ is a $k$-linear isomorphism $\mm/\mm^2 \to \mm/\mm^2$
such that $\Sym^2 \ell \colon \mm^2/\mm^3 \to \mm^2/\mm^3$
maps $V$ to $V'$
and $\eta_1,\ldots,\eta_d \in \frac{\mm^2}{\mm^3}/V'$.
(The isomorphism $A_V \to A_{V'}$ attached to such data is given by
$x_i \mapsto \ell(x_i) + \eta_i$.)
It follows that given $(V,\phi)$,
the dimension of the fiber of $Z_{d,r} \to \frakB_{n,k}$ containing it
equals the number of parameters needed to specify 
$\ell$ and $\eta_1,\ldots,\eta_d$.
Thus the fibers have dimension $d^2 + dr$.

Hence
\[
        \dim \frakB_{n,k} \ge \dim \pi(Z_{d,r}) 
        = r \left( \frac{d(d+1)}{2} - r \right) + n^2 - (d^2+dr).
\]
\end{proof}

\begin{theorem}
\label{T:lower bound}
For $n \ge 1$, 
\[
        \dim \frakB_{n,k} \ge 
        \begin{cases}
        \frac{2}{27} n^3 + \frac{1}{9} n^2 + \frac{5}{3} n - 1, &\text{if $n \equiv 0 \pmod{3}$} \\
        \frac{2}{27} n^3 + \frac{1}{9} n^2 + \frac{14}{9} n - \frac{20}{27}, &\text{if $n \equiv 1 \pmod{3}$} \\
        \frac{2}{27} n^3 + \frac{1}{9} n^2 + \frac{5}{3} n - \frac{37}{27}, &\text{if $n \equiv 2 \pmod{3}$.} \\
        \end{cases}
\]
\end{theorem}

\begin{proof}
Calculus shows that for fixed $n \ge 1$,
the bound in Lemma~\ref{L:Zdr} (with $r=n-1-d$)
as a function of $d$
increases up to a point around $2n/3-7/6+o(1)$
and decreases thereafter.
Evaluating the bound on the two integers on each side
shows that the maximum at nonnegative integers occurs at 
$d=\lfloor \frac{2n-2}{3} \rfloor$,
and the maximum value is as shown.
\end{proof}

\begin{remark}
The idea to use local algebras with $\mm^3=0$ to
obtain a lower bound of the form $(\frac{2}{27}+o(1))n^3$
is old, but Theorem~\ref{T:lower bound} seems to be more precise
than other results of its type.
\end{remark}

\begin{corollary}
\label{C:11}
If $n \ge 11$, then $\dim \frakB_{n,k} > n^2$.
\end{corollary}

\begin{remark}
It is natural to guess that 
the true dimension of $\frakB_{n,k}$ equals $n^2$ for $n \le 10$,
and equals the lower bound of Theorem~\ref{T:lower bound} for $n \ge 11$.
\end{remark}

\begin{proposition}
\label{P:8}
If $n \ge 8$, then $\frakB_{n,k}$ is reducible.
\end{proposition}

\begin{proof}
It suffices to give a rank-$8$ algebra that cannot be deformed to
an \'etale algebra, since then rank-$n$ algebras with the same property
can be obtained by taking a product with an \'etale algebra of rank $n-8$.
An explicit rank-$8$ algebra with this property was first given 
in \cite{Iarrobino-Emsalem1978}.
An example somewhat simpler than the one given there is
\[
	\frac{k[a,b,c,d]}{(a^2,ab,b^2,c^2,cd,d^2,ad-bc)},
\]
whose construction grew out of discussions of a working group at Berkeley
(Jonah Blasiak, Dustin Cartwright, David Eisenbud, Daniel Erman, 
Mark Haiman, the present author, Bernd Sturmfels, Mauricio Velasco,
and Bianca Viray).
\end{proof}

\begin{proposition}
\label{P:not order}
For each $n \ge 8$,
there exists a prime $p$ and an $\F_p$-algebra 
that is not isomorphic to $A/pA$
for any order $A$ in an \'etale $\Q$-algebra.
\end{proposition}

\begin{proof}
We will show that any algebra coming from a point of 
$\frakB_n(\F_p)-\overline{\frakB_n^\et}(\F_p)$ has the property;
the algebra
\[
	A_0:=\frac{\F_p[a,b,c,d]}{(a^2,ab,b^2,c^2,cd,d^2,ad-bc)} \times \F_p \times \cdots \times \F_p
\]
{}from the proof of Proposition~\ref{P:8},
with $n-8$ factors of $\F_p$, is an example.

If $A$ is an order in an \'etale $\Q$-algebra,
and we choose a $\Z$-basis,
then the corresponding section $s \colon \Spec \Z \to \frakB_n$
is such that $s(\Spec \Q) \in \frakB_n^\et$,
so $s(\Spec \Z) \subseteq \overline{\frakB_n^\et}$,
and $s(\Spec \F_p) \in \overline{\frakB_n^\et}(\F_p)$.
In other words, the point of $\frakB_n(\F_p)$ corresponding to $A/pA$
with its basis is in $\overline{\frakB_n^\et}(\F_p)$,
so $A/pA \not\isom A_0$.
\end{proof}

\begin{question}
Can one characterize the $k$-algebras that correspond
to points of $\overline{\frakB_{n,k}^\et}$?
\end{question}

Closely related is the following:
\begin{question}
What functor does the $\Z$-scheme $\overline{\frakB_n^\et}$ represent?
\end{question}

\section{Dimension: upper bound}
\label{S:dimension upper bound}

We now work toward an asymptotically matching upper bound 
on $\dim \frakB_{n,k}$,
namely $\dim \frakB_{n,k} \le \frac{2}{27}n^3 + O(n^{8/3})$.
Such a result is announced in \cite{Neretin1987}, 
who writes that the proof is nearly identical to the proof
he gives for the moduli space of $n$-dimensional Lie algebras.
We will give details of a proof for $\frakB_{n,k}$.

The approach towards both those results 
is to adapt the proof (begun in \cite{Higman1960}
and completed in \cite{Sims1965}) that the number of $p$-groups
of order $p^n$ is $p^{\frac{2}{27}n^3 + O(n^{8/3})}$.
As suggested to us by Hendrik Lenstra,
there is an analogy 
between the powers of the maximal ideal
of a local finite-rank $k$-algebra
and the descending $p$-central series of a $p$-group.
Although there seems to be no direct connection between 
finite-rank $k$-algebras and finite $p$-groups,
the combinatorial structure in the two enumeration proofs
are nearly identical.

\subsection{Symmetric bilinear maps}

This subsection is inspired by \cite{Sims1965}*{\S2}, which studied
alternating bilinear maps.

Throughout this subsection, we fix the following notation.
Let $V$ and $W$ be vector spaces over a field $k$.
Let $m=\dim V$ and $n=\dim W$.

\begin{proposition}
\label{P:good basis}
Let $(\;,\;)\colon V \times V \to W$ be a symmetric bilinear map.
Suppose $(V,V)=W$ but for no proper subspace $U \le V$ is $(U,U)=W$.
Then there exists a basis $x_1,\ldots,x_m$ of $V$
such that for $1 \le i \le m-1$, 
if $V_i$ is the span $\langle x_1,\ldots,x_i \rangle$
then $(x_i,x_{i+1}) \notin (V_i,V_i)$.
\end{proposition}

\begin{proof}
We show by induction that for $r=0,1,\ldots,m$ that there exist
$x_1,\ldots,x_r$ such that $(x_i,x_{i+1}) \notin (V_i,V_i)$
for $1 \le i \le r-1$.
If $r \le 1$, this is trivial.
The proof of the inductive step is as in the proof of 
\cite{Sims1965}*{Proposition~2.1}.
\end{proof}

\begin{corollary}
\label{C:m and n+1}
Under the hypothesis of Proposition~\ref{P:good basis}, $m \le n+1$.
\end{corollary}

\begin{corollary}
\label{C:small subspace generates}
For any $(\;,\;)$, there exists a subspace $U \le V$
of dimension $\le n+1$ with $(U,U)=(V,V)$.
\end{corollary}

\begin{proof}
Apply Corollary~\ref{C:m and n+1} to the minimal subspace $U$ satisfying
$(U,U)=(V,V)$.
\end{proof}

\begin{proposition}
\label{P:counting symmetric bilinear maps}
Suppose that $k$ is a finite field $\F_q$ of order $q$.
Let $V,W,m,n$ be as before,
and let $s \in \Z_{>0}$.
The number of symmetric bilinear maps $(\;,\;)\colon V \times V \to W$
such that for some $s$-dimensional subspace $U \le V$ we have $(U,U)=W$,
but for no $(s-1)$-dimensional subspace $U' \le V$ is $(U',U')=W$
is less than or equal to
\[
	q^{\frac{m^2}{2}(n-s) + O((m+n)^{8/3})}.
\]
\end{proposition}

\begin{proof}
See the proof of \cite{Sims1965}*{2.4--2.6}.
The only differences are
\begin{enumerate}
\item
We have $q$ instead of a prime $p$; this is of no consequence.
\item
Since we deal with symmetric forms in place of his alternating forms,
the $\binom{m}{2}$ that appears in the proof of Sims' Proposition~2.6 
should be replaced by $\binom{m+1}{2}$.
Both of these are close enough to $m^2/2$
that the difference in the resulting exponent
is dominated by the $O((m+n)^{8/3})$ term.
\end{enumerate}
\end{proof}

\subsection{Finite local $\F_q$-algebras}
\label{S:finite local algebras}

Our next step is to count rank-$n$ local $\F_q$-algebras $A$
with residue field $\F_q$.
Given such an $A$, we will define a collection of data,
and then prove that the data uniquely determine $A$ up to isomorphism.

Let $\mm$ be the maximal ideal of $A$.
For $i \ge 0$, let $m_i=\dim_{\F_q} \mm^i/\mm^{i+1}$.
Let $g_{01}=1$.
Let $V:=\mm/\mm^2$ and $W:=\mm^2/\mm^3$.
The multiplication in $A$ induces a symmetric bilinear map
$(\;,\;)\colon V \times V \to W$ such that $(V,V)=W$.
Let $\bar{V} \le V$ be a subspace of minimum dimension
such that $(\bar{V},\bar{V})=W$.
Let $s:=\dim \bar{V}$.
By Corollary~\ref{C:small subspace generates}, $s \le m_2+1$.
By Proposition~\ref{P:good basis}, we may choose a basis 
$x_1,\ldots,x_s$ of $\bar{V}$ such that 
for $i \in [1,s-1]$ we have $(x_i,x_{i+1}) \notin (V_i,V_i)$,
where $V_i:=\langle x_1,\ldots,x_i \rangle$.
For $i \in [0,s]$, let $W_i=(V_i,V_i)$ and $d_i=\dim W_i$.
So $0=d_0 \le d_1 < d_2 < \cdots < d_s = m_2$.
Choose a basis $y_1,\ldots,y_{m_2}$ of $W$ such that
$W_i=\langle y_1,\ldots,y_{d_i} \rangle$.
Given $i$ and $j \in (d_{i-1},d_i]$
we may assume $y_j=(x_b,x_i)$ 
for some $b \in [1,i]$ depending on $j$.
Extend $x_1,\ldots,x_s$ to a basis $x_1,\ldots,x_{m_1}$ of $V$.
Choose a representative $g_{1i} \in A$ of each $x_i$.
For each $i$ and each $j \in (d_{i-1},d_i]$,
take $g_{2j}=g_{1b} g_{1i}$ as the representative of $y_j$
for the same $b$ as above.
By induction on $i$, the map
\[
	\frac{\mm^{i-1}}{\mm^i} \tensor \bar{V} \to \frac{\mm^i}{\mm^{i+1}}
\]
induced by multiplication in $A$ 
is surjective for each $i \ge 2$.
Thus in particular, by induction on $i$,
for $i \ge 3$, we can choose $g_{i1},\ldots,g_{i,m_i} \in A$
representing a basis of $\mm^i/\mm^{i+1}$
such that each $g_{ij}$ equals $g_{i-1,r} g_{1\ell}$ 
for some $r \in [1,m_{i-1}]$ and $\ell \in [1,s]$.
For $i=3$, if $X_h$ denotes the image of 
\[
	\frac{\mm^2}{\mm^3} \tensor V_h \to \frac{\mm^3}{\mm^4}
\]
we may assume moreover that the $g_{3j}=g_{2r} g_{1\ell}$
are chosen so that the first few have $\ell=1$ and span $X_1$,
and the next few have $\ell=2$ and with the previous ones span $X_2$,
and so on up to $\ell=s$.
Then in each product $g_{2r} g_{1\ell}$ arising at the $\ell^\tH$ stage
(i.e., mapping to $X_\ell-X_{\ell-1}$), we have $r>d_{\ell-1}$,
since otherwise $g_{2r}=g_{1b} g_{1c}$ for some $b \le c \le \ell-1$
and the class of $g_{2r} g_{1\ell} = g_{1b} g_{1\ell} g_{1c}$
in $\mm^3/\mm^4$ lies in $X_{\ell-1}$, and hence is not a new basis
element of the $\ell^\tH$ stage.

The $g_{ij}$ for $i \ge 0$ and $j \in [1,m_i]$ form a basis for $A$.
Define $c_{i j \ell u v} \in \{0,1,\ldots,p-1\}$ by
\[
	g_{ij} g_{1\ell} = \sum c_{ij\ell uv} g_{uv}.
\]

\begin{proposition}
\label{P:local ring determined}
The isomorphism type of a local $\F_q$-algebra $A$ with residue field $\F_q$
is determined by the sequence $(m_i)_{i \ge 0}$
and the $c_{i j\ell u v}$ 
for $(i,j,\ell,u,v)$ satisfying the following conditions:
\begin{center}
\begin{tabular}{l}
$i=1$, $1 \le j \le \ell \le m_1$, $u \ge 2$, $1 \le v \le m_u$; or \\
$i=2$, $1 \le \ell \le s$, $d_{\ell-1} < j \le m_2$, $u \ge 3$, $1 \le v \le m_u$; or \\
$i\ge 3$, $1 \le j \le m_i$, $1 \le \ell \le s$, $u>i$, $1 \le v \le m_u$; \\
\end{tabular}
\end{center}
where $s$ and the $d_i$ are determined by the $c_{1j\ell u v}$.
\end{proposition}

\begin{proof}
Given the data above, 
we construct the vector space underlying $A$
by taking the $\F_q$-vector space with basis $g_{ij}$.
Since $A$ is generated as a $\F_q$-algebra by the $g_{1\ell}$ 
for $1 \le \ell \le m_1$,
it suffices to show that the given data
determine $g_{ij} g_{1\ell}$ for $i \ge 1$ 
as a linear combination of the $g_{uv}$.

If $i=1$, then by commutativity we may assume $j \le \ell$,
in which case the value of $g_{1j} g_{1\ell}$ is already given
by the $c_{1j\ell uv}$.

We next prove by strong induction on $\ell$ 
that $g_{2j} g_{1\ell}$ is determined for all $j$ and all $\ell \le s$.
Suppose the values $g_{2j} g_{1\ell'}$ 
have been determined for all $j$ whenever $\ell'<\ell$,
and we want to determine $g_{2j} g_{1\ell}$.
If $j>d_{\ell-1}$, the relevant $c_{2j\ell uv}$
have already been given.
So assume $j \le d_{\ell-1}$.
Then $g_{2j}=g_{1b} g_{1c}$ for some $b \le c \le \ell-1$,
so $g_{2j} g_{1\ell} = (g_{1b} g_{1\ell}) g_{1c}$.
Here $g_{1b} g_{1\ell}$ is a known combination of
the $g_{uv}$ for $u \ge 2$ and all $v$.
But $g_{2v} g_{1c}$ has been determined already
by the inductive hypothesis,
and $g_{uv} g_{1c}$ for $u \ge 3$ is already given
(since $c \le \ell-1 \le s$),
so the product $g_{2j} g_{1\ell} = (g_{1b} g_{1\ell}) g_{1c}$
is determined.

We now know that multiplication by $g_{1\ell}$ on all the $g_{ij}$ 
with $i \ge 1$ is determined provided that $\ell \le s$.
To extend this to all $\ell \le m_1$, we use induction on $i$.
The case $i=1$ is already given by the known $c_{1j\ell uv}$.
So assume $i \ge 2$.
Then $g_{ij}=g_{i-1,r} g_{1b}$ for some $r,b$ with $b \le s$,
so $g_{ij} g_{1\ell} = (g_{i-1,r} g_{1\ell}) g_{1b}$.
By the inductive hypothesis, $g_{i-1,r} g_{1\ell}$
is a known combination of $g_{uv}$ with $u \ge 2$,
so its product with $g_{1b}$ is determined (since $b \le s$).
\end{proof}

\begin{proposition}
\label{P:local k-algebras with residue field k}
The number of rank-$n$ local $\F_q$-algebras with residue field $\F_q$
is $q^{\frac{2}{27} n^3 + O(n^{8/3})}$
as $n \to \infty$.
The implied constant can be chosen independent of $q$.
\end{proposition}

\begin{proof}
The lower bound comes from Lemma~\ref{L:Zdr} and Theorem~\ref{T:lower bound}:
using standard formulas for $\#\Gr(r,d(d+1)/2)(\F_q)$
and $\#\GL_n(\F_q)$, one obtains a lower bound
that is at least a constant times $q$ raised to the right hand side
in Theorem~\ref{T:lower bound}.

Now for the upper bound.
To each local $\F_q$-algebra $A$ with residue field $\F_q$ we associate
the finite list of positive integers $m_0,m_1,\ldots,m_t$ summing to $n$,
and the $c_{ij\ell uv} \in \F_q$
as in Proposition~\ref{P:local ring determined}.
The number of possibilities for $m_0,m_1,\ldots,m_t$,
including the choice of $t$,
is at most $2^{n-1}$ (given $n$ stars, place or do not place a bar between
each consecutive pair of stars);
this is $O(q^{n^{8/3}})$, and hence may be ignored.

Below, ``log'' means logarithm to the base $q$.
The number of possibilities for $s$ is at most $n+1$,
so we may assume $s$ is fixed.
For each $s$, 
Proposition~\ref{P:counting symmetric bilinear maps}
shows that the log of the 
number of possibilities for the $c_{1j\ell 2v}$ giving
that value of $s$ is 
\[
	\frac{m_1^2}{2}(m_2-s) + O((m_1+m_2)^{8/3})
	= \frac{m_1^2}{2}(m_2-s) + O(n^{8/3}).
\]
The log of the number of possibilities for the $c_{1j\ell uv}$ with $u \ge 3$
is $\frac{m_1(m_1+1)}{2} (n-1-m_1-m_2)$.
Thus the log of the number of possibilities for all the $c_{1j\ell uv}$ is
\[
	\frac{m_1^2}{2}(n-1-m_1-s) + O(n^{8/3}).
\]

The log of the number of possibilities for the $c_{2j\ell uv}$ is
\[
	\left( \sum_{\ell=1}^s (m_2-d_{\ell-1}) \right) (n-1-m_1 - m_2).
\]
Since $0 \le d_0 \le d_1<d_2<\ldots<d_s$, 
we have $d_{\ell-1} \ge \ell-2$,
so the value above is
\[
	\le \left( m_2 s - \frac{s^2}{2} \right) (n-1-m_1 - m_2) + O(n^2).
\]

The log of the number of possibilities for the $c_{ij\ell uv}$ 
for a {\em fixed} $i \ge 3$ is
\[
	m_i s (n-1-m_1 - \cdots -m_i).
\]

Thus the log of the number of possibilities for {\em all} the $c_{ij\ell uv}$ 
is
\[
	\le \frac{m_1^2}{2}(n-1-m_1-s) + \left( m_2 s - \frac{s^2}{2} \right) (n-1-m_1 - m_2) + \sum_{i \ge 3} m_i s (n-1-m_1 - \cdots -m_i) + O(n^{8/3}),
\]
which up to $O(n^2)$ is the same as the expression $M$ 
in \cite{Sims1965}*{p.~165}.
Thus by the analysis there it is at most $\frac{2}{27} n^3 + O(n^{8/3})$.
\end{proof}

\begin{proposition}
\label{P:local Fq-algebras}
The number of rank-$n$ local $\F_q$-algebras
is $q^{\frac{2}{27} n^3 + O(n^{8/3})}$
as $n \to \infty$.
The implied constant can be chosen independent of $q$.
\end{proposition}

\begin{proof}
The lower bound follows from 
Proposition~\ref{P:local k-algebras with residue field k}.
Each local $\F_q$-algebra has residue field $\F_{q^d}$ for some $d \mid n$,
so by Proposition~\ref{P:local k-algebras with residue field k}
their total number is
\[
	\sum_{d \mid n} (q^d)^{\frac{2}{27} (n/d)^3 + O((n/d)^{8/3})}.
\]
This sum has at most $n$ terms, and each is at most 
$q^{\frac{2}{27} n^3 + O(n^{8/3})}$,
so the result follows.
\end{proof}

\subsection{Algebras of finite rank}
\label{S:algebras of finite rank}

\begin{lemma}
\label{L:convex}
Suppose $f\colon \R_{\ge 0} \to \R$ is convex and $f(0)=0$.
Then $f(x_1+\cdots+x_n) \ge f(x_1)+\cdots+f(x_n)$
for any $x_1,\ldots,x_n \in \R_{\ge 0}$.
\end{lemma}

\begin{proof}
This is a special case of the 
Hardy-Littlewood-P\'olya majorization 
inequality \cite{Hardy-Littlewood-Polya1988}*{Theorem~108}, 
for instance.
\end{proof}

\begin{theorem}
\label{T:Fq-algebras}
The number of rank-$n$ $\F_q$-algebras
is $q^{\frac{2}{27} n^3 + O(n^{8/3})}$
as $n \to \infty$.
The implied constant can be chosen independent of $q$.
\end{theorem}

\begin{proof}
The lower bound follows from Proposition~\ref{P:local Fq-algebras}.
A general $\F_q$-algebra is a product of local ones.
Thus we may specify a rank-$n$ $\F_q$-algebra ring by giving a partition 
$n = \lambda_1 + \cdots + \lambda_m$
and a local $\F_q$-algebra $A_i$ of rank $\lambda_i$ for each $i$.

The number of partitions (even if we did not impose an ordering on the
$\lambda_i$) is $\le 2^{n-1} < q^{O(n^{8/3})}$,
so it will suffice to bound the number of algebras for a fixed partition.
By Proposition~\ref{P:local Fq-algebras},
this number is
\[
	\le q^{f(\lambda_1)+\cdots+f(\lambda_m)}
\]
where $f(x):=\frac{2}{27} x^3 + c x^{8/3}$
for some universal constant $c>0$.
Lemma~\ref{L:convex} gives the desired upper bound.
\end{proof}

\begin{theorem}
\label{T:upper bound}
We have $\dim \frakB_{n,k} = \frac{2}{27} n^3 + O(n^{8/3})$
uniformly in $k$.
\end{theorem}

\begin{proof}
We may replace $k$ by its minimal subfield.
For fixed $n$, $\dim \frakB_{n,\Q} = \dim \frakB_{n,\F_p}$ 
for all but finitely many primes $p$,
so it suffices to obtain a bound for $\dim \frakB_{n,\F_p}$
that is uniform in $p$.
Theorem~\ref{T:Fq-algebras} estimates $\#\frakB_{n,\F_p}(\F_{p^e})$
for every $e$ (the latter count includes the choice of a basis,
but the number of choices is only $O(p^{n^2})$).
Taking $e \to \infty$ and applying the Lang-Weil bounds \cite{Lang-Weil1954},
we obtain the desired result.
\end{proof}

\section{Commutative rings of finite order}
\label{S:finite rings}

As a bonus, Theorem~\ref{T:Fq-algebras} leads to an asymptotic formula for
the number of (commutative) rings of order $p^n$,
namely $p^{\frac{2}{27}n^3 + O(n^{8/3})}$.
To prove this, we follow \cite{Kruse-Price1970},
which proved the analogous formula
$p^{\frac{4}{27}n^3 + O(n^{8/3})}$
for the number of associative rings of order $p^n$
that are not necessarily commutative or unital.
We begin with the commutative analogue of 
\cite{Kruse-Price1970}*{Theorem~3.1}.

\begin{lemma}
\label{L:reduction}
The number of (commutative) rings of order $p^n$ up to isomorphism
is at most $p^{n^2+n}$ times the number of rank-$n$ $\F_p$-algebras 
up to isomorphism.
\end{lemma}

\begin{proof}
(\cite{Kruse-Price1970})
For each ring $R$ of order $p^n$, 
choose generators $x_1,\ldots,x_m$ of the additive group of $R$
such that their orders $p^{a_1},\ldots,p^{a_m}$ multiply to $p^n$.
For $i \le m$ and $0\le j < a_i$, let $y_{ij}=p^j x_i$.
Rename all the $y_{ij}$ in any order as $z_1,\ldots,z_n$.
Then $z_i z_j = \sum c_{ijk} z_k$ for some $c_{ijk} \in \{0,1,\ldots,p-1\}$.
Construct the rank-$n$ $\F_p$-algebra $A$ having the same
structure constants $c_{ijk}$ considered in $\F_p$.
Associativity, commutativity, and existence of $1$ for $R$
imply the corresponding properties for $A$.

The construction above defines 
a map from the set of isomorphism classes of rings of order $p^n$ 
to the set of pairs $(\vec{a},A)$ where $\vec{a}=(a_1,\ldots,a_m)$
is a sequence of positive integers summing to $n$,
and $A$ is a based rank-$n$ $\F_p$-algebra.
Reversing the construction shows that this map is injective.

The number of sequences $\vec{a}$ is $2^{n-1} \le p^n$,
and the number of choices of basis for a rank-$n$ $\F_p$-algebra
is $\#\GL_n(\F_p) \le p^{n^2}$.
\end{proof}

\begin{theorem}
\label{T:rings of order p^n}
The number of commutative rings of order $p^n$ up to isomorphism
is $p^{\frac{2}{27} n^3 + O(n^{8/3})}$
as $n \to \infty$.
The implied constant can be chosen independent of $p$.
\end{theorem}

\begin{proof}
Combine Theorem~\ref{T:Fq-algebras} and Lemma~\ref{L:reduction}.
\end{proof}

\begin{theorem}
\label{T:rings of order up to N}
The number of commutative rings of order $\le N$ up to isomorphism
is $\exp(\frac{2}{27} (\log N)^3 / (\log 2)^2 + O((\log N)^{8/3}))$
as $N \to \infty$.
\end{theorem}

\begin{proof}
It suffices to prove the same bound for the number $r_N$ of rings 
of {\em exact} order $N$,
since $r_0+\cdots+r_N \le (N+1) \max_{i \in [0,N]} r_i$,
and $N+1 < \exp((\log N)^{8/3})$ for large $N$.

By Theorem~\ref{T:rings of order p^n}, there exists $c \ge 0$ such that
the number of rings of order $p^n$
is $\le p^{(2/27)n^3 + c n^{8/3}}$.
Write $N=\prod_p N_p$, where $N_p$ is a power of the prime $p$.
Then
\begin{align*}
	\log r_N 
	&\le \sum_{p \mid N}  \left( \frac{2}{27} \left(\frac{\log N_p}{\log p} \right)^3 + c \left(\frac{\log N_p}{\log p} \right)^{8/3} \right) \log p \\
	&\le \sum_{p \mid N}  \left( \frac{2}{27} \frac{(\log N_p)^3}{(\log 2)^2} + c \frac{(\log N_p)^{8/3}}{(\log 2)^{5/3}} \right) \\
	&\le \frac{2}{27} \frac{(\log N)^3}{(\log 2)^2} + c \frac{(\log N)^{8/3}}{(\log 2)^{5/3}},
\end{align*}
by Lemma~\ref{L:convex}.
\end{proof}

The main contribution to the number of finite rings
may come from those of the form
$(\frac{\Z_2[[x_1,\ldots,x_d]]}{\mm^3})/V$
where $\Z_2$ is the ring of $2$-adic integers, 
$\mm=(2,x_1,\ldots,x_d)$, and $V$ is an $\F_2$-subspace of $\mm^2/\mm^3$.
If this is so, then we would obtain the following information
about the ``typical'' finite ring:
\begin{conjecture}
As $N \to \infty$, the fraction of (isomorphism classes of) 
commutative rings $A$ of order $\le N$
satisfying the following conditions tends to $1$:
\begin{enumerate}
\item
The size of $A$ is $2^n$, where $2^n$ is the largest power of $2$ less than
or equal to $N$.
\item
The ring $A$ is local.
\item
The residue field of $A$ is $\F_2$.
\item
If $\mm$ is the maximal ideal of $A$,
then $\mm^3=0$ but $\mm^2\ne 0$.
\item
The $\F_2$-dimension of $\mm/\mm^2$ is $2n/3 + O(1)$.
\item
The characteristic of $A$ is $8$ (i.e., we have $8=0$ but $4 \ne 0$ in $A$).
\end{enumerate}
\end{conjecture}

Similar questions have been raised for groups:
see \cite{Mann1999}, for instance.

\section{Hilbert schemes revisited}
\label{S:Hilbert schemes revisited}

Finally, we use the results and techniques in 
Sections \ref{S:dimension lower bound} and~\ref{S:dimension upper bound}
to estimate the dimension of $\Hilb_n(\Aff^d)$ when $d$ is at least
about as large as $n$.

\begin{corollary}
\label{C:Hilbert scheme dimension 2}
If $d \ge n-1$,
\[
	\dim \Hilb_n(\Aff^d)_k = \frac{2}{27}n^3 + nd + O(n^{8/3}).
\]
\end{corollary}

\begin{proof}
Substitute Theorem~\ref{T:upper bound}
in Corollary~\ref{C:Hilbert scheme dimension}.
\end{proof}

We may generalize Corollary~\ref{C:Hilbert scheme dimension 2}
to cover a larger range of $d$.

\begin{theorem}
\label{T:Hilbert scheme dimension in larger range}
Suppose $d=\alpha n$ for some $\alpha \in \Q$.
Set
\[
	c_\alpha := \begin{cases} 
	   2/27, &\text{ if $\alpha \ge 2/3$} \\
	   \frac{\alpha^2}{2}(1-\alpha) &\text{ if $0 \le \alpha \le 2/3$.}
	\end{cases}
\]
Then for $\alpha \ge 1/2$,
\[
	\dim \Hilb_n(\Aff^d)_k = c_\alpha n^3 + nd + O(n^{8/3}).
\]
For $\alpha<1/2$, we have only the inequality
\[
	\dim \Hilb_n(\Aff^d)_k \ge c_\alpha n^3 + O(n^2).
\]
\end{theorem}

\begin{proof}
Let $\frakB_{n,k,d}$ be the image of $\frakH_n(\Aff^d) \to \frakB_n$,
so $\frakB_{n,k,d}$ parameterizes rank-$n$ algebras that can be
generated by $d$ elements.
By Corollary~\ref{C:open locus in Bn}, $\frakB_{n,k,d}$ is open
in $\frakB_{n,k}$.
The proof of Corollary~\ref{C:Hilbert scheme dimension}
generalizes to show that
\[
	\dim \Hilb_n(\Aff^d)_k = \dim \frakB_{n,k,d} - n^2 + nd,
\]
so it suffices to estimate $\dim \frakB_{n,k,d}$.

For $d \ge 2n/3$, we still have 
\[
	\dim \frakB_{n,k,d} = \frac{2}{27} n^3 + O(n^{8/3}),
\]
since the algebras constructed in Lemma~\ref{L:Zdr} 
for the lower bound on $\dim \frakB_{n,k}$ in Theorem~\ref{T:lower bound}
could be generated by about $2n/3$ elements.

For $d = \alpha n$ with $\alpha<2/3$,
we can use Lemma~\ref{L:Zdr} to obtain
\[
	\dim \frakB_{n,k,d} \ge c_\alpha n^3 - O(n^2).
\]
If $1/2 \le \alpha \le 2/3$, we obtain a matching upper bound 
(but with error term $O(n^{8/3})$) by generalizing the proof of
Theorem~\ref{T:upper bound} as follows.
Proceed as in the proof of 
Proposition~\ref{P:local k-algebras with residue field k},
but at the end impose the additional condition that $m_1 \le \alpha n$.
This translates, in the notation of \cite{Sims1965}*{p.~165},
into the constraint $x \le \alpha$ in addition to 
the constraints $x+y \le 1$, $0\le y \le x$.
The constant in front of the $n^3$ in the lower bound obtained is
the maximum value of
\[
	B(x,y) := x^2(1-x-y)/2 + y^2(1-x-y)/2 + y(1-x-y)^2/2
\]
subject to these constraints.
For fixed $x \ge 1/2$, calculus shows that the maximum 
value of $B(x,y)$ for $y \in [0,1-x]$ 
occurs at $y=0$, and the value there is $\frac{x^2}{2}(1-x)=c_x$.
Calculus shows also that 
the maximum value of $B(x,y)$ on the triangle with vertices
$(0,0)$, $(1/2,0)$, $(1/2,1/2)$ is attained at $(1/2,0)$.
By the previous two sentences, for $1/2 \le \alpha \le 2/3$,
the maximum value of $B(x,y)$ on the region defined by
$x \le \alpha$ in addition to 
the constraints $x+y \le 1$, $0\le y \le x$
is $\frac{\alpha^2}{2}(1-\alpha)$, attained at $(\alpha,0)$.
\end{proof}

\begin{remark}
Before the present work, it seems that asymptotic bounds
for $\dim \Hilb_n(\Aff^d)_k$ were known only as $n \to \infty$ for fixed $d$:
it is shown in \cite{Briancon-Iarrobino1978} that
this dimension is bounded above and below by universal positive constants
times $n^{2-2/d}$.
\end{remark}

\section*{Acknowledgements} 

I thank Hendrik Lenstra for the insight that
enumeration of finite-rank algebras 
might be related to enumeration of $p$-groups.
I thank also Manjul Bhargava, Bas Edixhoven, Mark Haiman, and David Helm 
for discussions.
This work was begun at the trimester on ``Explicit methods in number theory''
at the Institut Henri Poincar\'e in 2004,
and continued at the workshop on ``Rings of low rank''
at the Lorentz Center in 2006,
in addition to U.\ C.\ Berkeley,
supported partially by NSF grant DMS-0301280.

\end{document}